\documentclass[reqno,11pt]{amsart}

\usepackage{amsmath}
\usepackage{amsfonts}
\usepackage{amssymb}
\usepackage{color}
\usepackage{mathrsfs}
\usepackage[colorlinks=true]{hyperref}

\usepackage{cleveref}

\newtheorem{thm}{Theorem}[section]

\newtheorem{cor}[thm]{Corollary}

\numberwithin{equation}{section}

\newcommand{\R}{\mathbb{R}}

\newcommand{\N}{\mathbb{N}}

\newcommand{\A}{\mathbb{A}}
\newcommand{\B}{\mathbb{B}}

\newcommand{\E}{\mathbb{E}}

\newcommand{\ml}{\mathcal{L}}

\newcommand{\ve}{\varepsilon}
\newcommand{\rd}{\mathrm{d}}
\newcommand{\dom}{\mathrm{dom}}

\newcommand{\bear}{\begin{eqnarray}} 
\newcommand{\eear}{\end{eqnarray}} 
\newcommand{\bean}{\begin{eqnarray*}} 
\newcommand{\eean}{\end{eqnarray*}} 
\newcommand{\bs}{\begin{split}}
\newcommand{\es}{\end{split}}

\newcommand{\dhr}{\mathrel{\lhook\joinrel\relbar\kern-.8ex\joinrel\lhook\joinrel\rightarrow}}

\begin{document}

\title[Compact Resolvent and Perturbations]{A Note on the Compactness of the Resolvent of the Age-Diffusion Operator}

\author{Christoph Walker}
\email{walker@ifam.uni-hannover.de}
\address{Leibniz Universit\"at Hannover\\ Institut f\" ur Angewandte Mathematik \\ Welfengarten 1 \\ D--30167 Hannover\\ Germany}
 
\date{\today}

\begin{abstract}
The generator of the semigroup associated with linear age-structured population models including spatial diffusion is shown to have compact resolvent. 
\end{abstract}

\keywords{Age structure, diffusion, semigroups.}
\subjclass[2010]{47D06, 35B35, 35M10, 92D25}

\maketitle

\section{Introduction and Main Results}\label{App:CompRes}

The dynamics of a population consisting of individuals structured by age and spatial position is governed by equations of the form
\begin{subequations}\label{PPPP} 
\begin{align}
\partial_t u+ \partial_au \, &=     A(a)u \,, \qquad t>0\, ,\quad a\in (0,a_m)\, ,\label{P1}\\ 
u(t,0)&=\int_0^{a_m}b(a)\, u(t,a)\, \rd a\,, \qquad t>0\, ,\label{P2} \\
u(0,a)&=  u_0(a)\,, \qquad a\in (0,a_m)\,,\label{P3} 
\end{align}
\end{subequations}
where  $u=u(t,a)$ denotes the population density 
 $
u:\R^+\times J\rightarrow E_0
$
with values in a Banach space $E_0$ representing the spatial heterogeneity.
The age of individuals is denoted by $a\in J:=[0,a_m]$ with
 maximal age $a_m\in (0,\infty)$. Spatial movement of individuals is described by the age-dependent (usually: differential) operator~$A$, for which we assume that
 there is $\rho>0$ with	
\begin{subequations}\label{A}
\begin{equation}\label{A1}
A\in  C^\rho\big(J,\mathcal{H}(E_1,E_0)\big)\,.
\end{equation}
$\mathcal{H}(E_1,E_0)$ stands for the space of generators of analytic semigroups on $E_0$ with domain~$E_1$ equipped with the operator norm in $\ml(E_1,E_0)$.
 Death processes of individuals are incorporated into the operator $A$ and spatial boundary conditions into the domain of definition~$E_1$.  We impose that 
\begin{equation}\label{A0}
E_1\stackrel{d}{\dhr} E_0\,,
\end{equation}
that is, $E_1$ is a densely and compactly embedded subspace of the Banach space $E_0$.
Equation~\eqref{P2} is the birth law with birth rate $b$, for which we shall assume that
\begin{equation}\label{A2l}
b\in L_{\infty}\big(J,\ml(E_\theta)\big)\,, \quad \theta\in  \{0,\vartheta\}  
\,,
\end{equation} 
\end{subequations}
for some fixed $\vartheta\in (0,1)$. Here, the space $E_\vartheta:= (E_0,E_1)_\vartheta$ is an interpolation space with admissible interpolation functor $(\cdot,\cdot)_\vartheta$ (see~\cite{LQPP}). In fact, given any admissible interpolation functor $(\cdot,\cdot)_\theta$ with $\theta\in (0,1)$ we use the notion $E_\theta:= (E_0,E_1)_\theta$ and set $$\E_\theta:=L_1(J,E_\theta)\,,\quad \theta\in [0,1]\,.$$
In~\eqref{P3}, $u_0$ represents the initial population. The main features of~\eqref{PPPP}  are that the evolution equation~\eqref{P1} has a hyperbolic character if spatial movement is neglected (i.e. $A=0$) or a parabolic character if aging is neglected and that condition~\eqref{P2} is nonlocal  with respect to age.

Age-structured population models (with and without spatial diffusion) have been the subject of intensive research in the past;  we refer  e.g. to \cite{WebbSpringer} and the references therein  for more details and concrete examples.

Realistic models for the evolution of age- and spatially structured populations are nonlinear and involve density-dependent vital rates $b$ or nonlinear  operators~$A$ incorporating nonlinear death processes or nonlinear diffusion. The study of linear equations like~\eqref{PPPP} is important for the understanding of such more complex models. For instance, they arise as (part of the) linearization when investigating stability properties of equilibria~\cite{Walker_Stability2023}. It is well-known that the solutions to the linear problem~\eqref{PPPP} may be represented by a semigroup on the phase space \mbox{$\E_0=L_1(J,E_0)$}. We refer to \cite{WebbSpringer,Thieme_DCDS98,Rhandi98,RhandiSchnaubelt_DCDS99,WalkerIUMJ,WalkerJEPE} and the references therein for the investigation of this semigoup  though this list is far from being complete.\\ 

It follows from \cite[Theorem~2.3~(c)]{WalkerIUMJ} that the generator $\A:\dom(\A)\subset \E_0\to \E_0$ of this semigroup is given by
\begin{align*}
\dom(\A):=\bigg\{\psi\in C(J,E_0)\,;\,&  \exists\, \zeta_\psi\in \E_0\, \text{such that $\psi$ is the mild solution to}\\
&\partial_a\psi =  A(a) \psi -\zeta_\psi(a)\,,\  a\in J\,,\ \ \psi(0)=\int_0^{a_m} b(a) \psi(a)\,\rd a
 \bigg\}
\end{align*}
and
$$
\A \psi := \zeta_\psi\,.
$$
Mild solution here means that
$$
\psi(a)=\Pi(a,0)\psi(0)-\int_0^a\Pi(a,\sigma)\,\zeta(\sigma)\,\rd \sigma\,,\quad a\in J\,,
$$
where  
$$
\big\{\Pi(a,\sigma)\in\ml(E_0)\,;\, a\in J\,,\, 0\le\sigma\le a\big\} 
$$ 
denotes the parabolic evolution operator on $E_0$ with regularity subspace $E_1$ in the sense of \cite[Section~II.2.1]{LQPP} (existence and uniqueness of $\Pi$ is ensured by~\eqref{A1} and~\cite[II.~Corollary~4.4.2]{LQPP}). 
With this notation, the linear problem~\eqref{PPPP} can be reformulated as the Cauchy problem
\begin{equation}\label{CP}
u'(t)=\A u(t) \,,\quad t\ge 0\,,\qquad u(0)=u_0
\end{equation}
in $\E_0$. Since $\A$ generates a strongly continuous semigroup $(e^{t\A})_{t\ge 0}$ on~$\E_0$, there is for each initial value $u_0\in \dom(\A)$ a unique solution 
$$
u\in C^1(\R^+,\E_0)\cap C(\R^+,D(\A))
$$
to~\eqref{CP}, where $D(\A)$ means the domain $\dom(\A)$ equipped with the graph norm. Precise information on $D(\A)$ -- which is defined implicitly initially -- is thus pertinent, also with regard to evolution problems involving time-dependent operators $\A=\A(t)$, see \cite{WalkerDCDS23}. For instance, one can show \cite[Lemma~2.1]{WalkerDCDS23} that  
$$
\left\{\psi\in \E_1\cap W_1^1(J,E_0)\,;\,\psi(0)=\int_0^{a_m}b(a)\, \psi(a)\, \rd a\right\}
$$
is a core for $\A$; that is, a dense subspace of $D(\A)$. Owing to~\eqref{A0} and \cite[Corollary~4]{Simon_CompactSets}, this space is also compactly embedded into $\E_0$. Moreover, from \cite[Theorem~2.3~(c)]{WalkerIUMJ} it follows  that  $D(\A)\hookrightarrow \E_\theta$ for each~$\theta\in [0,1)$, where the spaces $\E_\theta$ play an important role in the study of nonlinear problems as domains of definitions for the nonlinearities.
Herein we prove that  these embeddings are compact:


\begin{thm}\label{CompResolv}
Suppose \eqref{A}. Then  $\A$ has compact resolvent. In fact, for every~$\theta\in [0,1)$, the embedding $D(\A)\hookrightarrow  \E_\theta$  is compact. 
\end{thm}


Compact resolvent means that $(\lambda-\A)^{-1}$ is a compact operator on $\E_0$ for every $\lambda>0$ large enough. In order to prove Theorem~\ref{CompResolv}, the definition of $D(\A)$ entails to derive compactness of the mapping
$$
\zeta\mapsto \int_0^a\Pi(a,\sigma)\,\zeta(\sigma)\,\rd \sigma
$$ 
 in $\E_0=L_1(J,E_0)$. To this end we shall follow the lines of~\cite{BarasHassanVeron}.\\

Besides compactness of the embedding $D(\A)\hookrightarrow \E_0$, Theorem~\ref{CompResolv} also yields information on the spectrum of $\A$, e.g. that it consists of eigenvalues only. This information was derived previously in \cite{WalkerIUMJ} by showing that the semigroup $(e^{t\A})_{t\ge 0}$ is eventually compact (under slightly stronger assumptions). However, for perturbations of the form $\A+\B$ (arising, for instance, when linearizing equations incorporating nonlinear vital rates), the eventual compactness of the corresponding semigroup has recently been obtained in \cite{Walker_Stability2023}, but only for particular perturbations~$\B$. Nonetheless, Theorem~\ref{CompResolv} ensures for general perturbations~$\B\in\ml (\E_\alpha,\E_0)$ with $\alpha\in [0,1)$ the compactness of the resolvent of $\A+\B$ and thus provides spectral properties for perturbed operators $\A+\B$:


\begin{cor}\label{Cor:PertCompResolv}
Suppose \eqref{A}. Let $\B\in\ml (\E_\alpha,\E_0)$ for some $\alpha\in [0,1)$. Then $\A+\B$ with domain \mbox{$\dom(\A)$} generates a strongly continuous semigroup $(e^{t(\A+\B)})_{t\ge 0}$ on~$\E_0$ and has compact resolvent. In particular, the spectrum $\sigma(\A+\B)$ is a pure point spectrum without finite accumulation point.
\end{cor}

If $E_0$ is an ordered Banach space and
\begin{subequations}\label{P}
\begin{equation}
 \text{$A(a)$ is resolvent positive for each $a\in J$}
\end{equation}
and 
\begin{equation}
b\in L_\infty(J,\ml_+(E_0))\,,
\end{equation}
then the semigroup $(e^{t\A})_{t\ge 0}$ generated by $\A$ is positive on $\E_0$, see \cite{WalkerIUMJ}; that is, $\A$ is resolvent positive.
\end{subequations}
In case that $E_0$ is a Banach lattice, more information on the spectral bound $$s(\A):=\sup\{\mathrm{Re}\,\lambda\,;\,\lambda\in\sigma(\A)\}$$ is available, see \cite{BFR}:


\begin{cor}\label{Cor2}
Let $E_0$ be a Banach lattice and suppose \eqref{A} and \eqref{P}. If $\B\in\ml_+ (\E_\alpha,\E_0)$ for some $\alpha\in [0,1)$, then $\A+\B$  generates a positive semigroup~$(e^{t(\A+\B)})_{t\ge 0}$ on $\E_0$ with $e^{t\A} \le e^{t(\A+\B)}$  for $t\ge 0$. Moreover, $s(\A)\le s(\A+\B)$ and 
$$
(\lambda-\A)^{-1}\le (\lambda-\A-\B)^{-1}\,,\quad\lambda>s(\A+\B)\,.
$$
If $s(\A+\B)>-\infty$, then $s(\A+\B)$ is an eigenvalue of $\A+\B$ possessing an eigenvector $\psi\in \dom(\A)$ with $\psi\ge 0$. E.g. this is the case if 
$b(a)\Pi(a,0)\in \ml_+(E_0)$ is strongly positive for $a$ in a subset of $J$ of positive measure.
\end{cor}



We prove the above statements in the next sections. To this end, we assume throughout the following~\eqref{A}.


\section{Proof of Theorem~\ref{CompResolv}}

\subsection*{I. Compact Resolvent}

First note that, for $\lambda\in\R$, the family
$$
\Pi_{\lambda}(a,\sigma):=e^{-\lambda (a-\sigma)}\Pi(a,\sigma)\,,\qquad a\in J\,,\quad 0\le\sigma\le a\,,
$$
is the evolution operator associated with $-\lambda+A$. Moreover,  according to \cite[II.~Lemma~5.1.3]{LQPP}  there are $\varpi\in\R$ and $M_\vartheta\ge 1$ with
	\begin{equation}\label{EO}
	\|\Pi(a,\sigma)\|_{\ml(E_\vartheta)}+(a-\sigma)^\vartheta\,\|\Pi(a,\sigma)\|_{\ml(E_0,E_\vartheta)}\le M_\vartheta e^{\varpi (a-\sigma)}\,,\quad 0\le \sigma\le a\le a_m \,.
	\end{equation}
We may choose $\lambda>0$ in the resolvent set of $\A$ so large such that for the operator
$$
Q_\lambda:=\int_0^{a_m} b(a)\Pi_{\lambda}(a,0)\,\rd a\in \ml(E_\vartheta)
$$
we have $\|Q_\lambda\|_{\ml(E_\vartheta)}<1$ (see~\eqref{A2l} and~\eqref{EO}).
Consider a bounded sequence $(\phi_j)_{j\in\N}$ in~$\E_0$. Then, for $j\in\N$,
$$
\psi_j:=(\lambda-\A)^{-1}\phi_j\in\mathrm{dom}(\A)\,, 
$$
and the definition of $\A$ implies
\begin{equation}\label{p10}
\psi_j(a) =\Pi_{\lambda}(a,0)\psi_j(0)+\int_0^a\Pi_{\lambda}(a,\sigma) \phi_j(\sigma)\,\rd \sigma\,,\quad a\in J\,,
\end{equation}
and
\begin{equation}\label{p20}
 (1-Q_\lambda)\psi_j(0)=  \int_0^{a_m} b(a)\int_0^a \Pi_{\lambda}(a,\sigma)\, \phi_j(\sigma) \,\rd \sigma\,\rd a \,.
\end{equation}
In particular, since by \eqref{A2l} and \eqref{EO} we have
\begin{align*}
\Bigg\|&\int_0^{a_m} b(a)\int_0^a \Pi_{\lambda}(a,\sigma)\, \phi_j(\sigma) \,\rd \sigma\,\rd a \Bigg\|_{E_\vartheta}
\\
&\le M_\vartheta \,\|b\|_{L_\infty(J,\ml(E_\vartheta))}\, \int_0^{a_m} \int_0^a e^{(-\lambda+\varpi) (a-\sigma)}(a-\sigma)^{-\vartheta}\, \|\phi_j(\sigma) \|_{E_0}\,\rd \sigma\,\rd a\\
&= M_\vartheta \,\|b\|_{L_\infty(J,\ml(E_\vartheta))}\,\int_0^{a_m} \|\phi_j(\sigma) \|_{E_0} \int_\sigma^{a_m} e^{(-\lambda+\varpi) (a-\sigma)} (a-\sigma)^{-\vartheta} \,\rd a\,\rd \sigma\le c \|\phi_j\|_{\E_0}\,,
\end{align*}
the sequence
$$
\left(\int_0^{a_m} b(a)\int_0^a \Pi_{\lambda}(a,\sigma)\, \phi_j(\sigma) \,\rd \sigma\,\rd a\right)_{j\in\N}
$$
is bounded in $E_\vartheta$. Thus, \eqref{p20} and $\|Q_\lambda\|_{\ml(E_\vartheta)}<1$ imply that 
\begin{equation}\label{pq}
(\psi_j(0))_{j\in\N}\ \text{  is bounded in}\ E_\vartheta\,.
\end{equation} 
We adopt the proof of~\cite{BarasHassanVeron} in order to prove that
the set $\{u_j\,;\, j\in\N\}$, given by
$$
u_j(a):=\int_0^a\Pi_{\lambda}(a,\sigma) \phi_j(\sigma)\,\rd \sigma\,,\quad a\in J\,,\quad j\in\N\,,
$$
is relatively compact in $\E_0$. We split this into two steps: \vspace{2mm}

{\bf (i)}  Let $\mu>0$ be fixed and define
$$
v_j^\mu(a):=\Pi_{\lambda}(\mu+a,a)u_j(a)=\int_0^a\Pi_{\lambda}(\mu+a,\sigma) \phi_j(\sigma)\,\rd \sigma\,,\quad a\in J\,,\quad j\in\N\,,
$$
where we used the evolution property\footnote{Set $A(a):=A(a_m)$ for $a> a_m$ to obtain $A\in C^\rho(\R^+,\mathcal{H}(E_1,E_0))$. The uniqueness assertion of~\cite[II.~Corollary 4.4.2]{LQPP} ensures that the corresponding evolution operator then also extends $\Pi$.}
$$
\Pi_{\lambda}(\mu+a,a)\Pi_{\lambda}(a,\sigma)=\Pi_{\lambda}(\mu+a,\sigma)\,,\quad \sigma\le a\le \mu+a \,.
$$
Analogously to the derivation of~\eqref{pq},  the sequence 
\begin{align}\label{G2}
\text{$(v_j^\mu)_{j\in\N}$ is bounded in $\E_\vartheta$}\,.
\end{align} 
Note that the space $E_\vartheta$ embeds compactly into $E_0$ due to~\eqref{A0}.

For the equi-integrability recall from  \cite[Section~II.2.1]{LQPP} that the evolution operator has the continuity property $\Pi\in C(\Delta_J^*,\ml(E_0))$, where $\Delta_J^*:=\{(a,\sigma)\in J\times J\,;\, \sigma < a\}$. Uniform continuity on compact subsets of $\Delta_J^*$ implies that, given $\ve>0$, there is $\eta:=\eta(\ve,\mu)>0$ such that
$$
\|\Pi_{\lambda}(a_1,\sigma_1)-\Pi_{\lambda}(a_2,\sigma_2)\|_{\ml(E_0)}\le \ve
$$
whenever $(a_i,\sigma_i)\in \Delta_J^*$, $\mu\le a_i-\sigma_i$ for $i=1,2$ and $\vert (a_1,\sigma_1)-(a_2,\sigma_2)\vert\le \eta$.
We use this and~\eqref{EO} to derive, for $h\in (0,\eta)$ with  \mbox{$0<a\le a+h\le a_m$},
\begin{align*}
\|v_j^\mu(a+h)-v_j^\mu(a)\|_{E_0}&\le
 \int_a^{a+h}\|\Pi_{\lambda}(\mu+a+h,\sigma)\|_{\ml(E_0)}\,\|\phi_j(\sigma)\|_{E_0}\,\rd\sigma\\
&\quad + \int_0^{a}\|\Pi_{\lambda}(\mu+a+h,\sigma)-\Pi_{\lambda}(\mu+a,\sigma)\|_{\ml(E_0)}\,\|\phi_j(\sigma)\|_{E_0}\,\rd\sigma\\
&\le c(a_m)\int_a^{a+h} \|\phi_j(\sigma)\|_{E_0}\,\rd\sigma+\ve  \int_0^{a}\|\phi_j(\sigma)\|_{E_0}\,\rd\sigma
\end{align*}
and therefore
\begin{align*}
\int_0^{a_m}\|\tilde v_j^\mu(a+h)-\tilde v_j^\mu(a)\|_{E_0}\,\rd a &\le
 h c(a_m)\|\phi_j \|_{\E_0}+\ve a_m \|\phi_j \|_{\E_0}\,,
\end{align*}
with  tilde referring to  trivial extension. Since $(\phi_j)_{j\in\N}$ is bounded in $\E_0$ and $\ve>0$ was arbitrary, we deduce
\begin{align*}
\lim_{h\to 0}\, \sup_{j\in\N}\int_0^{a_m}\|\tilde v_j^\mu(a+h)-\tilde v_j^\mu(a)\|_{E_0}\,\rd a =0\,.
\end{align*} 
That is, $\{v_j^\mu\,;\, j\in\N\}$  is equi-integrable. Therefore, taking into account~\eqref{G2} we are in a position to apply~\cite[Theorem~3]{Simon_CompactSets}  and derive that
\begin{align}\label{equi}
\{v_j^\mu\,;\, j\in\N\}\ \text{  is relatively compact in  $\E_0$ for $\mu>0$}\,.
\end{align} 

{\bf (ii)} We consider the limit $\mu\to 0$. Given arbitrary $\ve\in (0,a_m)$, we argue as in part {\bf (i)} to find $\eta(\ve)>0$  such that
$$
\|\Pi_{\lambda}(a_1,\sigma_1)-\Pi_{\lambda}(a_2,\sigma_2)\|_{\ml(E_0)}\le \ve
$$
whenever $(a_i,\sigma_i)\in \Delta_J^*$, $\ve\le a_i-\sigma_i$ for $i=1,2$ and $\vert (a_1,\sigma_1)-(a_2,\sigma_2)\vert\le \eta(\ve)$. Thus, for $0<\mu<\eta(\ve)$ and $a\in J$ with $a\ge \ve$, we obtain from~\eqref{EO} that
\begin{align*}
\|v_j^\mu(a)-u_j(a)\|_{E_0}&\le
 \int_0^{a-\ve}\|\Pi_{\lambda}(\mu+a,\sigma)-\Pi_{\lambda}(a,\sigma)\|_{\ml(E_0)}\,\|\phi_j(\sigma)\|_{E_0}\,\rd\sigma\\
&\quad +  \int_{a-\ve}^a \big(\|\Pi_{\lambda}(\mu+a,\sigma)\|_{\ml(E_0)}+\|\Pi_{\lambda}(a,\sigma)\|_{\ml(E_0)}\big)\,\|\phi_j(\sigma)\|_{E_0}\,\rd\sigma\\
&\le \ve \int_0^{a-\ve} \|\phi_j(\sigma)\|_{E_0}\,\rd\sigma + c(a_m)\int_{a-\ve}^a \|\phi_j(\sigma)\|_{E_0}\,\rd\sigma\,, 
\end{align*}
while, for $0\le a\le \ve$, we have
\begin{align*}
\|v_j^\mu(a)-u_j(a)\|_{E_0}&\le  c(a_m)\int_0^{\ve} \|\phi_j(\sigma)\|_{E_0}\,\rd\sigma\,.
\end{align*}
Therefore,
\begin{align*}
\|v_j^\mu-u_j\|_{\E_0}&\le  c(a_m)\, \ve \, \|\phi_j\|_{\E_0}+\ve\, a_m\, \|\phi_j\|_{\E_0}
\end{align*}
so that, since $\ve>0$ was arbitrary,
\begin{align*}
\lim_{\mu\to 0}\, \sup_{j\in\N}\, \|v_j^\mu-u_j\|_{\E_0}=0\,.
\end{align*}
Together with \eqref{equi} we conclude that $\{u_j\,;\, j\in\N\}$ is relatively compact in $\E_0$.\vspace{2mm}

{\bf (iii)} Finally,  since  
$$
\|\Pi_{\lambda} (a+h,0)-\Pi_{\lambda} (a,0)\|_{\ml(E_\vartheta,E_0)}\le c h^{\vartheta}\,,\quad 0\le a\le a+h\le a_m\,,
$$
according to \cite[II.Equation(5.3.8)]{LQPP}, we infer from \eqref{pq} and the Arzel\`a - Ascoli Theorem that $(\Pi_{\lambda} (\cdot,0)\psi_j(0))_{j\in\N}$ is relatively compact in $C(J,E_0)$. Consequently, the previous step~{\bf (ii)} and~\eqref{p10} entail that $(\psi_j)_{j\in\N}$ is relatively compact in $\E_0$. Therefore, $\A$ has compact resolvent.

\subsection*{II. Compact Embedding}

Since $\A$ has compact resolvent, the embedding $D(\A)\hookrightarrow \E_0$ is compact. The compact embedding into $\E_\theta$ for $\theta\in (0,1)$ fixed follows by an interpolation argument. Indeed, choosing $\ve>0$ with $\theta<1-3\ve$ we set $\theta_k:=\theta (1-k\ve)^{-1}$ for $k=1,2,3$ so that $0<\theta_1<\theta_2<\theta_3<1$. Then \cite[Equation I.(2.5.2)]{LQPP} yields for the continuous interpolation functor $(\cdot,\cdot)_{\theta_3,\infty}^0$ that
$$
\big(\E_0,\E_{1-\ve}\big)_{\theta_3,\infty}^0\hookrightarrow \big(\E_0,\E_{1-\ve}\big)_{\theta_2,1}=\big(L_1(J,E_0),L_1(J,E_{1-\ve})\big)_{\theta_2,1}
$$
while \cite[Theorem 1.18.4]{Triebel} ensures for the real interpolation functor $(\cdot,\cdot)_{\theta_2,1}$ that
$$
\big(L_1(J,E_0),L_1(J,E_{1-\ve})\big)_{\theta_2,1}\doteq L_1\big(J,(E_0,E_{1-\ve})_{\theta_2,1}\big)\,.
$$
Finally, \cite[I.~Remark~2.11.2~(a)]{LQPP} implies
$$
L_1\big(J,(E_0,E_{1-\ve})_{\theta_2,1}\big)\hookrightarrow L_1\big(J,E_{(1-\ve)\theta_1}\big)=\E_\theta\,.
$$
Gathering these findings we obtain the continuous embedding
$$
\big(\E_0,\E_{1-\ve}\big)_{\theta_3,\infty}^0\hookrightarrow \E_\theta
$$
so that there is $c>0$ with
\begin{equation}\label{w}
\|\phi\|_{\E_\theta}\le c \|\phi\|_{\E_0}^{1-\theta_3}\, \|\phi\|_{\E_{1-\ve}}^{\theta_3}\,,\quad \phi\in \E_{1-\ve}\,. 
\end{equation}
Now, if $(\phi_j)_{j\in\N}$ is a bounded sequence in $D(\A)$, then it is also bounded in $\E_{1-\ve}$ due to the continuous embedding \mbox{$D(\A)\hookrightarrow \E_{1-\ve}$}. Moreover, since the embedding \mbox{$D(\A)\hookrightarrow \E_0$} is compact, we may assume without loss of generality that it is a Cauchy sequence in~$\E_0$. According to~\eqref{w}  it is also a Cauchy sequence in $\E_\theta$ and thus converges. Consequently, the embedding $D(\A)\hookrightarrow \E_\theta$ is compact. This proves Theorem~\ref{CompResolv}.\qed


\section{Proof of Corollary~\ref{Cor:PertCompResolv}}


It was shown in \cite[Theorem 2.1]{WalkerIUMJ} that $\A+\B$ with domain \mbox{$\dom(\A)$} generates a strongly continuous semigroup on $\E_0$ when $\B\in\ml (\E_\alpha,\E_0)$ for some $\alpha\in [0,1)$.
Thus, for $\lambda>0$ sufficiently large, both $\lambda-\A$ and $\lambda-\A-\B$ are invertible with
\begin{equation*}
(\lambda-\A-\B)^{-1}= (\lambda-\A)^{-1}\big(1-\B(\lambda-\A)^{-1}\big)^{-1}\,.
\end{equation*}
Since $(\lambda-\A)^{-1}\in\ml(\E_0)$ is compact according to Theorem~\ref{CompResolv} while $$\big(1-\B(\lambda-\A)^{-1}\big)^{-1}\in\ml(\E_0)$$ since $D(\A) \hookrightarrow \E_\alpha$ (see \cite[Theorem 2.3]{WalkerIUMJ}), it follows that $(\lambda-\A-\B)^{-1}$ is compact.   
This proves Corollary~\ref{Cor:PertCompResolv}.\qed


\section{Proof of Corollary~\ref{Cor2}}


Let $E_0$ be a Banach lattice and assume, in addition to \eqref{A}, also~\eqref{P}.  Corollary~\ref{Cor2} follows from~\cite[Proposition~12.11]{BFR} (except that therein the perturbation $\B\in\ml_+(\E_0)$ is a bounded operator on $\E_0$), the proof is the same. Indeed, it was shown for $\B\in\ml_+(\E_\alpha,\E_0)$ in \cite[Theorem~1.2]{WalkerIUMJ} that the semigroup $(e^{t(\A+\B)})_{t\ge 0}$ generated by $\A+\B$ is positive and satisfies
\begin{equation*}\label{E33T}
e^{t(\A+\B)}\phi=e^{t\A}\phi+\int_0^te^{(t-s)\A)}\,\B\, e^{s(\A+\B)}\,\phi\,\rd s\,,\quad t\ge 0\,,\quad \phi\in\E_0\,.
\end{equation*}
Therefore, $e^{t\A} \le e^{t(\A+\B)}$  for $t\ge 0$. We then argue as in \cite[Proposition~12.11]{BFR} that the Laplace transform formula~\cite[Theorem~12.7]{BFR} implies
$$
(\lambda-\A)^{-1}\le (\lambda-\A-\B)^{-1}\,,\quad\lambda>\max\{s(\A+\B),s(\A)\}\,,
$$
and then $s(\A+\B)\ge s(\A)$ using \cite[Corollary~12.9]{BFR}.  

Finally, that $s(\A+\B)>-\infty$ is an eigenvalue with positive eigenvector follows from the Krein-Rutman Theorem~\cite[Theorem~12.15]{BFR} since $\A+\B$ is resolvent positive and has a compact resolvent according to Corollary~\ref{Cor:PertCompResolv}. In particular, if 
$b(a)\Pi(a,0)\in \ml_+(E_0)$ is strongly positive for $a$ in a subset of $J$ of positive measure, then $s(\A)\in \R$ according to~\cite[Proposition~4.2]{WalkerIUMJ}. 
This proves Corollary~\ref{Cor2}.\qed




\bibliographystyle{siam}
\bibliography{AgeDiff_230419}

\end{document}